\documentclass[10pt]{amsart}
\usepackage{amsmath, amssymb}
 \usepackage{mathrsfs}
\newcommand{\no}[1]{#1}
\renewcommand{\no}[1]{}
\no{\usepackage{times}\usepackage[subscriptcorrection, slantedGreek, nofontinfo]{mtpro}
\renewcommand{\Delta}{\upDelta}}
\usepackage{color}

\newtheorem{theorem}{Theorem}

\title[Lipschitz stability]{Corrigendum \vskip .3 cm An inverse problem in corrosion detection:\\ 
stability estimates, J. Inv. Ill-posed Problems\\ 12 (4) (2004), 349-367.
}

\thanks{The author is very grateful to Daijun Jiang and Jun Zou for their valuable comments during his stay at the Chinese University of Hong Kong on February 2017. He warmly thank the Chinese University of Hong Kong for hospitality.}

\author[Mourad Choulli]{Mourad Choulli}
\address{IECL, UMR CNRS 7502, Universit\'e de Lorraine, Boulevard des Aiguillettes BP 70239 54506 Vandoeuvre Les Nancy cedex- Ile du Saulcy - 57 045 Metz Cedex 01 France}
\email{mourad.choulli@univ-lorraine.fr}

\date{\today}

\begin{document}

\maketitle


Unless otherwise stated, $\Omega$ is a $C^\infty$ bounded domain of $\mathbb{R}^2$ so that its boundary $\Gamma $ is the union of two disjoint closed subsets with nonempty interior, $\Gamma =\Gamma_1 \cup \Gamma_2$.

\smallskip
We considered in \cite{Ch1} the stability issue  for the problem of determining the boundary coefficient $q$, appearing in the BVP
\begin{equation}\label{bvp}
\left\{
\begin{array}{ll}
\Delta u=\;\; &\mbox{in}\; \Omega ,
\\
\partial _\nu u +q u =0  &\mbox{on}\; \Gamma_1 ,
\\
\partial _\nu u =f &\mbox{on}\; \Gamma_2,
\end{array}
\right.
\end{equation}
from the boundary measurement $u_{|\gamma _2}$, where $\gamma_2$ is an open subset of $\Gamma_2$.

\smallskip
Our proof of \cite[Theorem 2.1]{Ch1} is partially incorrect. We rectify here this proof. We precisely establish a stability estimate of logarithmic type for the inverse problem described above. Contrary to the result announced in \cite[Theorem 2.1]{Ch1}, we do not know whether Lipschitz stability, even around a particular unknown coefficient, is true. Note that Lipschitz stability around an arbitrary unknown boundary coefficient is false in general as shows the following counter example in which $\Omega =\{1/2<|x|<1\}$, $\Gamma _1=\{|x|=1/2\}$ and $\Gamma _2=\{|x|=1\}$. Let, in polar coordinates system $(r,\theta )$,
\begin{align*}
&u=1+\ln r
\\
&u_k=u+2^{-k}k^{-2}(r^k+r^{-k})\cos (k\theta ),\;\; k\ge 1.
\end{align*}

By straightforward computations we check that $u$ and $u_k$ are the solutions of the BVP \eqref{bvp} respectively when
\begin{align*}
&q=\frac{2}{1-\ln 2},
\\
&q=q_k=\frac{2+k^{-1}(2^{-2k+1}-2)\sin (k\theta )}{1-\ln 2+k^{-2}(2^{-2n}+1)\sin (k\theta )},\;\; k\ge 1,
\end{align*}
and $f=1$.

\smallskip
By simple calculations, we get $\| u-u_k\|_{L^2(\Gamma _2)}=O\left(2^{-k}k^{-2}\right)$,  while $\| q-q_k\|_{L^2(\Gamma _1)}=O(k^{-1})$.

\smallskip
To our knowledge, the only case where Lipschitz stability holds is when $q$ is assumed to be a priori piecewise constant. We refer to \cite{Si} for more details.

\smallskip
Throughout, the unit ball of a Banach space $X$ is denoted by $B_X$ and 
\[
L^p_{K}(D)=\{ h\in L^p(D);\; \mbox{supp}(h)\subset K\},\;\; 1\le p\le \infty .
\]

For sake of clarity, we start our analysis with stability around a particular boundary coefficient. To this end, fix $0<\alpha <1$ and, for $0\le f\in C^{1,\alpha}(\Gamma _2)$, denote by $w(f)\in C^{2,\alpha}(\overline{\Omega})$ the solution of the BVP
\[
\left\{
\begin{array}{ll}
\Delta w=0\;\; &\mbox{in}\; \Omega ,
\\
w=0  &\mbox{on}\; \Gamma_1 ,
\\
\partial _\nu w=f &\mbox{on}\; \Gamma_2.
\end{array}
\right.
\]
According to the strong maximum principle and Hopf's lemma (see for instance \cite{GT}), $\partial _\nu w<0$ on $\Gamma_1$.


\smallskip
Let $q_0 =- \partial _\nu w(f)_{|\Gamma _1} (>0)$ and set $u_0=1+w$. Then it is straightforward to check that $u_0$ is the unique solution of the BVP
\[
\left\{
\begin{array}{ll}
\Delta u=0\;\; &\mbox{in}\; \Omega ,
\\
\partial _\nu u +q_0 u =0  &\mbox{on}\; \Gamma_1 ,
\\
\partial _\nu u =f &\mbox{on}\; \Gamma_2.
\end{array}
\right.
\]

For $(\varphi _1,\varphi_2) \in L^2(\Gamma _1)\oplus L^2(\Gamma _2)$, define $L (\varphi _1,\varphi _2) :=y$, where $y\in H^{3/2}(\Omega )$ is the unique weak solution of the BVP
\[
\left\{
\begin{array}{ll}
\Delta y=0\;\; &\mbox{in}\; \Omega ,
\\
\partial _\nu y +q_0 y =\varphi _1  &\mbox{on}\; \Gamma_1 ,
\\
\partial _\nu y =\varphi _2  &\mbox{on}\; \Gamma_2.
\end{array}
\right.
\]
An application of Green's formula leads 
\begin{align}
\int_\Omega |\nabla y|^2dx + \int_{\Gamma_1} q_0y ^2d\sigma &=\int_{\Gamma _1}\varphi_1 y d\sigma +\int_{\Gamma _2}\varphi_2 y d\sigma \label{0.0}
\\
&\le \|(\varphi_1,\varphi_2) \|_{L^2(\Gamma _1)\oplus L^2(\Gamma _2)}\|y \|_{H^1(\Omega )}.
\end{align}
Using that
\[
h\rightarrow \left(\int_\Omega |\nabla h|^2dx +\int_{\Gamma_1}q_0h ^2d\sigma \right)^{1/2}
\]
defines an equivalent norm on $H^1(\Omega )$, we derive from \eqref{0.0}
\begin{equation}\label{0.1}
\|y \|_{H^1(\Omega )}\le \kappa_0 \| (\varphi_1,\varphi_2)\|_{L^2(\Gamma _1)\oplus L^2(\Gamma _2)},
\end{equation}
for some constant $\kappa _0$ depending only on $\Omega$ and $f$.

\smallskip
As $y$ is also the solution of the BVP
\[
\left\{
\begin{array}{ll}
\Delta y=0\;\; &\mbox{in}\; \Omega ,
\\
\partial _\nu y + y = (1-q_0) y  +\varphi_1 &\mbox{on}\; \Gamma_1 ,
\\
\partial _\nu y =\varphi_2  &\mbox{on}\; \Gamma_2,
\end{array}
\right.
\]
we get from the usual a priori estimates for non homogenous BVP's (see \cite{LM}) that there exits a constant $\kappa _1$, depending only on $\Omega$ and $f$, so that
\[
\|y \|_{H^{3/2}(\Omega )}\le \kappa _1\|\ (\varphi_1,\varphi_2)\|_{L^2(\Gamma _1)\oplus L^2(\Gamma _2)}.
\]
In other words, we proved that $L \in \mathscr{B}(L^2(\Gamma _2),H^{3/2}(\Omega ))$ and
\begin{equation}\label{0.2}
\|L\|:=\| L \|_{\mathscr{B}(L^2(\Gamma _1)\oplus L^2(\Gamma _2),H^{3/2}(\Omega ))}\le \kappa _1.
\end{equation}

For $q\in L^2(\Gamma _1)$, define the operator $H_q$ as follows
\[
H_q:H^{3/2}(\Omega )\rightarrow H^{3/2}(\Omega ): H_q(u)=L \left(-qu_{|\Gamma _1},0\right).
\]
If $\kappa$ is the norm of the trace operator 
\[
h\in H^{3/2}(\Omega )\rightarrow u_{|\Gamma _1}\in C(\Gamma _1), 
\]
then
\[
\|H_q\|_{\mathscr{B}(H^{3/2}(\Omega ))}\le \kappa  \|L\|\|q\|_{L^2(\Gamma _1)}.
\]
Whence, for any $q\in \mathcal{U} =(2\kappa \|L\|)^{-1}B_{L^2(\Gamma _1)}$, $I-H_q$ is invertible and 
\begin{equation}\label{0.3}
\|\left(I-H_q\right)^{-1}\|_{\mathscr{B}(H^{3/2}(\Omega ))}\le 2,\;\;  q\in \mathcal{U}.
\end{equation}

Define, for $q\in \mathcal{U}$ and $(\varphi _1,\varphi_2) \in L^2(\Gamma _1)\oplus L^2(\Gamma _2)$, 
\[
u_q(\varphi _1,\varphi_2)=\left(I-H_q\right)^{-1}L (\varphi _1,\varphi_2).
\]

In light of the identity
\[
u_q(\varphi _1,\varphi_2)=L \left(-qu_{|\Gamma _1}+\varphi_1 ,\varphi_2\right),
\]
we derive  that $u_q(\varphi _1,\varphi_2)\in H^{3/2}(\Omega )$ is the solution of the BVP
\[
\left\{
\begin{array}{ll}
\Delta u=0\;\; &\mbox{in}\; \Omega ,
\\
\partial _\nu u + (q_0 +q)u = \varphi_1   &\mbox{on}\; \Gamma_1 ,
\\
\partial _\nu u =\varphi_2  &\mbox{on}\; \Gamma_2.
\end{array}
\right.
\]
Note that according to \eqref{0.3}
\begin{equation}\label{0.4}
\|u_q(\varphi _1,\varphi_2)\|_{H^{3/2}(\Omega )}\le 2 \kappa_1\| (\varphi_1,\varphi_2)\|_{L^2(\Gamma _1)\oplus L^2(\Gamma _2)}.
\end{equation}
\smallskip
Set $u_q=u_q(0,f)$. That is $u_q$ is the solution of the BVP
\[
\left\{
\begin{array}{ll}
\Delta u=0\;\; &\mbox{in}\; \Omega ,
\\
\partial _\nu u + (q_0 +q)u = 0   &\mbox{on}\; \Gamma_1 ,
\\
\partial _\nu u =f  &\mbox{on}\; \Gamma_2.
\end{array}
\right.
\]
Observe that  \eqref{0.4} yields
\begin{equation}\label{0.4.0}
\|u_q\|_{H^{3/2}(\Omega )}\le 2 \kappa_1\|f\|_{L^2(\Gamma _2)}.
\end{equation}

Let $\gamma _1$ be a nonempty open subset of $\Gamma_1$ so that $\Gamma _1\setminus \overline{\gamma_1}$ is nonempty. Define $L_{\overline{\gamma_1}} ^2(\Gamma_1)$ as the set of those functions $p\in L^2(\Gamma )$ so that $\mbox{supp}(p)\subset \overline{\gamma _1}$. We can mimic the proof of \cite[Propoistion 2.1]{Ch1} to show that the mapping
\[
\Phi :q\in \mathcal{U} \cap L_{\overline{\gamma_1}} ^2(\Gamma_1)\rightarrow \chi_{\Gamma _1}\left[\partial _\nu u_q{_{|\gamma _1}}\right]\in L_{\overline{\gamma_1}} ^2(\Gamma_1)
\]
is continuously Fr\'echet differentiable and $\Phi '(0)=N$. Here, for $p\in L_{\overline{\gamma_1}}^2(\Gamma _1)$, $Np=\chi_{\Gamma _1}\left[ \partial _\nu v_p{_{|\gamma _1}}\right]$, where $v_p$ is the solution of the BVP
\[
\left\{
\begin{array}{ll}
\Delta v=0\;\; &\mbox{in}\; \Omega ,
\\
\partial _\nu v + q_0 v= -p   &\mbox{on}\; \Gamma_1 ,
\\
\partial _\nu v =0  &\mbox{on}\; \Gamma_2.
\end{array}
\right.
\]

Similarly to the proof of \cite[Lemma 2.1]{Ch1}, we prove that $N$ is an isomorphism. Therefore, by the implicit function theorem,  there exists $\widetilde{\mathcal{U}} \subset \mathcal{U}$ so that $\Phi ^{-1} $ is Lipschitz continuous, on $\widetilde{\mathcal{V}}= \Phi  (\widetilde{\mathcal{U}} \cap L_{\overline{\gamma_1}} ^2(\Gamma_1))$, with Lipschitz constant less or equal to $2\|N^{-1}\|$. That is
\begin{equation}\label{5}
\|q_1-q_2\|_{L^2(\Gamma _1 )}\le 2\|N^{-1}\| \|\partial _\nu u_{q_1}-\partial_\nu u_{q_2 }\|_{L^2(\gamma _1)},\;\; q_1,q_2\in \widetilde{U} \cap L_{\overline{\gamma_1}} ^2(\Gamma_1) .
\end{equation}

\smallskip
Let $k$ be a positive 
integer, $s\in \mathbb{R}$, $1\leq r\leq \infty$ and consider the vector space
\[
B_{s,r}(\mathbb{R}^k ):=\{ w\in \mathscr{S}'(\mathbb{R}^k);\; (1+|\xi |^2)^{s/2}\widehat{w}\in 
L^r(\mathbb{R}^k )\},
\]
where $\mathscr{S}'(\mathbb{R}^k)$ is the space of temperated distributions on $\mathbb{R}^k$ and
$\widehat{w}$ is the Fourier transform of $w$. Equipped with the norm 
\[
\|w\|_{B_{s,r}(\mathbb{R}^k )}:=\left\|(1+|\xi |^2)^{s/2}\widehat{w}\right\|_{ L^r(\mathbb{R}^k )},
\]
$B_{s,r}(\mathbb{R}^k )$ is a Banach space. Note that $B_{s,2}(\mathbb{R}^k )$ is merely the Sobolev
space $H^s(\mathbb{R}^k )$. Using local charts and a partition of unity, we construct 
$B_{s,r}(\Gamma _1)$ from $B_{s,r}(\mathbb{R})$ similarly as $H^s(\Gamma _1)$ is built 
from $H^s(\mathbb{R})$.

\smallskip
Fix $m>0$. If $f\in H^{3/2}(\Gamma _2)$ and $q\in mB_{B_{3/2,1}(\Gamma _1)}$, then by \cite[Theorem 2.3]{Ch0}, $u_q\in H^3(\Omega )$ and
\begin{equation}\label{6}
\|u_q\|_{H^3(\Omega )}\le C_0.
\end{equation}
Here and henceforth, $C_0$ is a constant depending only on $\Omega$, $f$ and $m$.

\smallskip
But in dimension two $H^3(\Omega )$ is continuously embedded in $C^2(\overline{\Omega})$. Whence, \eqref{6} entails
\begin{equation}\label{7}
\|u_q\|_{C^2(\overline{\Omega})}\le C_0.
\end{equation}

Let
\[
\Psi (\rho )=|\ln \rho |^{-1/2}+\rho,\;\; \rho >0,
\]
extended by continuity at $0$ by setting $\Psi (0)=0$.

\smallskip
Let $\gamma_2$ be a nonempty open subset of $\Gamma _2$. According to \cite[Proposition 2.7]{Ch3}, there exists a constant $C>0$, depending only on $\Omega$, $f$, $m$ and $\gamma_2$, so that
\begin{equation}\label{8}
\|\partial _\nu u_{q_1}-\partial_\nu u_{q_2 }\|_{L^2(\gamma _1)} \le C\Psi \left( \|u_{q_1}-u_{q_2}\|_{H^1(\gamma _2)} \right).
\end{equation}

Set
\[
\mathscr{Q}_m=mB_{B_{3/2,1}(\Gamma _1)}\cap \widetilde{\mathcal{U}} \cap L_{\overline{\gamma_1}} ^2(\Gamma_1).
\]
Note that $\mathscr{Q}_m\neq \emptyset$ if $m$ is chosen sufficiently large. 

\smallskip
We can now combine \eqref{5} and \eqref{8} in order to obtain
\[
\|q_1-q_2\|_{L^2(\Gamma _1 )}\le C\Psi \left( \|u_{q_1}-u_{q_2}\|_{H^1(\gamma _2)} \right),\;\; q_1,q_2\in \mathscr{Q}_m.
\]

We sum up our analysis in the following theorem, where we used the fact that $H^{3/2}(\Gamma _2)$ is continuously embedded in $C^2(\Gamma _2)$,

\begin{theorem}\label{theorem1}
Let $0\le f\in H^{3/2}(\Gamma _2)$, $q_0 =- \partial _\nu w(f)_{|\Gamma _1} (>0)$ and $\gamma_i $ be  a nonempty open subset of $\Gamma_i$, $i=1,2$, with $\Gamma \setminus \overline{\gamma_1}\not= \emptyset$. There exists  a neighborhood $\widetilde{\mathcal{U}}$ of $q_0$ in $L^2(\Gamma _1)$, depending on $f$, $\Omega$ and $\gamma_1$ with the property that, if $m>0$ is chosen in such a way that
\[
\mathscr{Q}_m=mB_{B_{3/2,1}(\Gamma _1)}\cap \widetilde{\mathcal{U}} \cap L_{\overline{\gamma_1}} ^2(\Gamma_1)\neq \emptyset,
\]
we find a constant $C>0$, depending on $f$, $m$, $\Omega$ and $\gamma_i$, $i=1,2$, so that
\[
\|q_1-q_2\|_{L^2(\Gamma _1 )}\le C\Psi \left( \|u_{q_1}-u_{q_2}\|_{H^1(\gamma _2)} \right),\;\; q_1,q_2\in \mathscr{Q}_m.
\]
\end{theorem}

\smallskip
We now discuss briefly the stability around an arbitrary $q_0$. Let then $q_0\in L^\infty (\Gamma _1)$ be non negative and non identically equal to zero and let $f\in L^2(\Gamma _2)$ be non identically equal to zero. Denote by $u_0\in H^{3/2}(\Omega )$ the solution of the BVP
\[
\left\{
\begin{array}{ll}
\Delta u=0\;\; &\mbox{in}\; \Omega ,
\\
\partial _\nu u + q_0u = 0   &\mbox{on}\; \Gamma_1 ,
\\
\partial _\nu u =f  &\mbox{on}\; \Gamma_2.
\end{array}
\right.
\]
As it is observed in \cite{Ch1}, 
\[
\Gamma _0=\{x\in \Gamma _1;\; u_0(x)\ne 0\}
\]
is an open dense subset of $\Gamma_1$. 

\smallskip 
Slight modifications of the preceding analysis allow us to prove the following result 
\begin{theorem}\label{theorem2}
Let $f\in H^{3/2}(\Gamma _2)$, $f\not\equiv 0$, $0\le q_0\in L^\infty (\Gamma _1)$, $q_0\not \equiv 0$, $K$ a compact subset of $\Gamma _0$ so that $\Gamma _1\setminus K\not=\emptyset$ and $\gamma_2$ be  a nonempty open subset of $\Gamma_2$. There exists  a neighborhood $\widetilde{\mathcal{U}}$ of $q_0$ in $L^2(\Gamma _1)$, depending on $f$, $\Omega$ and $K$ with the property that, if $m>0$ is chosen in such a way that
\[
\mathscr{Q}_m=mB_{B_{3/2,1}(\Gamma _1)}\cap \widetilde{\mathcal{U}} \cap L_K^2(\Gamma_1)\neq \emptyset ,
\]
we find a constant $C>0$, depending on $f$, $m$, $\Omega$, $K$ and $\gamma_2$, so that
\[
\|q_1-q_2\|_{L^2(\Gamma _1 )}\le C\Psi \left( \|u_{q_1}-u_{q_2}\|_{H^1(\gamma _2)} \right),\;\; q_1,q_2\in \mathscr{Q}_m.
\]
\end{theorem}

Observe that, as in \cite{Ch1}, the last theorem can be extended to the case where $\partial \Gamma_1\cap \partial \Gamma _2\neq\emptyset$.

\smallskip
In the most general case, in dimensions two and three, we can prove a stability estimate of triple logarithmic type (see \cite[Theorem 4.9]{Ch3}).

\vskip .5cm
\end{document}